\documentclass[12pt]{article}
\usepackage{amssymb}
\usepackage{amsthm}
\usepackage{amsmath}
\usepackage{graphicx}
\graphicspath{{./Figs/}}
\usepackage{lmodern}
\usepackage{graphics}
\usepackage{float}

\newtheorem{teo}{Theorem}
\newtheorem{lem}{Lemma}

\newcommand{\ds}{\displaystyle}
\newcommand{\R}{\mathbb R}
\newcommand{\wm}{\mathrm{w}}
\newcommand{\pa}{\partial}

\newcommand{\vp}{\varphi}
\newcommand{\ve}{\varepsilon}
\newcommand{\om}{\omega}
\newcommand{\be}{\begin{equation}}
\newcommand{\ee}{\end{equation}}

\newcommand{\ipd}{\stackrel{\normalfont\text{def}}{=}}

\newcommand{\ap}{\alpha}
\newcommand{\bt}{\beta}
\newcommand{\gm}{\gamma}

\newcommand{\E}{\text{Er}}
\newcommand{\D}{\Delta}

\providecommand{\pare}[1]{\left(#1\right)}

\providecommand{\corc}[1]{\left[#1\right]}
\usepackage[usenames]{color}
\begin{document}
\allowdisplaybreaks

\title{A finite difference scheme for smooth solutions of the general mKdV equation}
\author{J.~Noyola Rodriguez
\thanks{Universidad Autónoma de Guerrero, Carlos E. Adame 54, 39650 Acapulco, Guerrero. Mexico,\
20264@uagro.mx}\and  G.~Omel'yanov\thanks{Corresponding author,
Universidad de Sonora, Rosales y Encinas, 83000 Hermosillo, Sonora, Mexico,\ omel@mat.uson.mx}}

\date{}
\maketitle
\begin{abstract}
We consider a generalization of the mKdV  model of shallow water out-flows. This generalization is a family of equations with nonlinear dispersion terms containing, in particular, KdV, mKdV, Benjamin-Bona-Mahony, Camassa-Holm, and Degasperis-Procesi equations.  Nonlinear dispersion, generally speaking, implies instability of classical solutions and wave breaking in a finite time. However, there are  special conditions under which  the general mKdV equation admits classical solutions that are global in time. We have  created an economic  finite difference scheme that  preserves this property for numerical solutions. To illustrate this we demonstrate some numerical results about propagation and interaction of solitons.
\end{abstract}

\emph{Key words}:  general mKdV equation,  Degasperis-Procesi model, finite difference scheme, soliton, interaction

\emph{2010 Mathematics Subject Classification}: 35Q35, 35Q53, 65M06

\section{Introduction}
 We consider a modern unidirectional approximation of the shallow water system  (see e.g. \cite{DegProc}) called the  ``general mKdV equation" (gmKdV):
\begin{align}
&\frac{\pa }{\pa t}\left\{u-\alpha^2\ve^2\frac{\pa^2 u}{\pa x^2}\right\}\label{1}\\
&+\frac{\pa}{\pa x}\left\{c_0u+c_1u^n
-c_2\ve^2\Big(\frac{\pa u}{\pa x}\Big)^2+\ve^2\big(\gamma-c_3u\big)\frac{\pa^2 u}{\pa x^2}\right\}=0, \; x \in \mathbb{R}^1, \; t > 0.\notag
\end{align}
Here $\alpha$, $c_0,\dots,c_3$, $\gamma$ are  real parameters  and  $\varepsilon$ characterizes the dispersion, $n=2$ or $n=3$.  The constants $\alpha\geq0$ and $\gamma\geq0$ are associated with different characters of the dispersion manifestation. In the Green-Naghdi approximation  the restriction $\alpha+\gamma>0$ is required. The equation (\ref{1}) terms with $c_2\geq0$ and $c_3\geq0$ can be treated as representations of nonlinear dispersion.  In the Camassa-Holm approximation  $c_2+c_3>0$.

This six parametric family of third order conservation laws contains as particular cases a list of basic equations: the KdV and mKdV equation if $\alpha=c_2=c_3=0$ and $n=2$ or $n=3$; the  Benjamin-Bona-Mahony (BBM) equation (\cite{BBM}, 1972) if $n=2$,
$\gamma=c_2=c_3=0$; the Camassa-Holm (CH) equation (\cite{CH}, 1993) if $n=2$, $c_2=c_3/2$, $c_1=3c_3/2\alpha^2$, and $\gamma=0$; and the Degasperis-Procesi (DP) equation (\cite{DegProc}, see also \cite{ConLan}) if $n=2$, $c_2=c_3$, $c_1=2c_3/\alpha^2$, and $c_0=\gamma=0$. All these particular equations are quite different. Indeed, the KdV, BBM and generalizations of CH and DP (if $\gamma+\alpha^2c_0>0$) equations  have soliton-type traveling wave solutions. At the same time,  the CH with $c_0=0$ and DP equations, under the condition $u\to0$ as $x\to\pm\infty$, have non-smooth traveling wave solutions only. Next, the KdV, mKdV, CH, and DP equations are completely integrable, whereas all others particular cases of the model (\ref{1}) are essentially non-integrable (see e.g. \cite{ELY}). Consequently,  KdV, mKdV, DP, and CH solitons collide elastically, whereas BBM ``solitons" have changed after the interaction and an oscillatory tail is generated \cite{BPS}. Furthermore, the Cauchy and periodic problems for the CH and DP equations have been studied extensively (see e.g. \cite{ConLan, ELY, Mus, Wahl} and references therein), whereas the solvability of similar problems for the general case (\ref{1}) remains be unknown.
So, the general model (\ref{1}) represents the non trivial  object of investigations and it is naturally to expect a very interesting behavior of its solution.

Concerning the numerical modeling of the equations of family (\ref{1}), it seems that quite a lot of researches has been carried out. The main part of them is devoted to the CH equation.
The fact is that when $c_3=2c_2$, the main term that generates instability falls out of the main balance law
\be
\frac{d }{d t}\left\{\int_{-\infty}^\infty u^2dx+\alpha^2\int_{-\infty}^\infty(\ve u_x)^2dx\right\}
=\ve^{-1}(c_3-2c_2)\int_{-\infty}^\infty(\ve u_x)^3dx.\label{4}
\ee
Further, due to wave breaking and the existence of non-smooth solutions, the proposed schemes for CH and DP equations are mainly based on Fourier representation or Galerkin method (e.g. \cite{KalLen, MatYam, Mat, LHY} and others). Obviously, this involves the use of fairly complex iterative procedures.

As for the finite difference approximation, the approach that is mainly developed here is one that uses bi-Hamiltonian structures of CH and DP equations  - the so-called ``mean vector field method" (\cite{CelGrim, MM, MMF, MF, Coclit} and others, see also \cite{Feng}).

Our interest in the numerical investigation of  the gDP model (for $n=2$) appeared after the first result \cite{Om2, OmNo} which states that gDP solitons,  under some conditions,  interact elastically in the weak asymptotic (for $\ve\to0$) sense. To demonstrate this effect numerically there has been developed an approach \cite{GO1, GO2, OmNo1} which adapts to the finite-difference representation the conservation law
\be
\frac{d}{d t}\int_{-\infty}^\infty u(x,t)\,dx=0,\label{3}
\ee
and the balance law (\ref{4}). Next, let us note that right-hand side in (\ref{4}) disappears for even functions. Thus, numerical simulation  of soliton motion seems to be stable.
Subsequent numerical experiments showed the validity of this hypothesis.
Moreover, it turned out that the collision of solitons, in which, generally speaking, the evenness of the solution is violated, does not lead to any significant errors.

In this paper we adapt the main ideas of \cite{OmNo1} and previous articles \cite{GO1, GO2}
 to the general mKdV model (\ref{1}) for $n=3$ and create a "conservative" and effective algorithm.

The content of the paper is the following: in Section 2 we present  assumptions which guarantee the soliton solution existence,  a description of the finite difference scheme for sufficiently  smooth waves is contained in Section 3, which focuses on the problem of dynamics and interaction of solitons. Section 4 shows the results of the corresponding numerical experiments.
In the appendix, we demonstrate some technical details of the finite difference scheme analysis.
\section{Soliton type solution}

The difference between gDP and gmKdV equations is much deeper than between KdV and mKdV. The most important  novelty here is such that solitons and antisolitons have different shapes. Furthermore, both solitons and antisolitons move with positive velocities.

To clarify the method for constructing a smooth traveling wave for equation (\ref{1}), we briefly explain the approach \cite{OmNo2}. Let's consider the ansatz
\begin{equation}\label{15}
u=A\om\big(\beta(x-Vt-x^0)/\varepsilon,A\big),
\end{equation}
which we call soliton (for $A>0$ or antisoliton for $A<0$), regardless of the scenario of their interaction. Here $\om(\eta,A)$ is a smooth function such that
\begin{equation}\label{15a}
\om(0,A)=1,\quad \om(-\eta,A)=\om(\eta,A), \quad \om(\eta,A)\to 0\quad \text{as}\quad \eta\to\pm\infty,
\end{equation}
 the wave amplitude $A\neq0$, the scale $\beta$, and the initial point $x^0$ are free parameters. The velocity $V=V(A)\neq0$ should be determined.

We assume
\be
\gm+\ap^2V>0, \label{7}
\ee
and, to simplify  subsequent formulas, define the notation
\begin{align}
&\begin{aligned}
&r=c_3/(c_2+c_3), \; \bt=\sqrt{c_1(\gm+\ap^2V)}/(c_3\sqrt{r}),\\
&q=c_3^2(V-c_0)/\big(c_1(\gm+\ap^2V)^2\big),\; p=c_3A/(\gamma+\alpha^2V).\label{8}
\end{aligned}
\end{align}
Next we define a new function $g=g(\eta,q)$ such that
\begin{equation}\label{12}
\om(\eta,A)=\big(1-g(\eta,q)^r\big)/p.
\end{equation}
Substituting (\ref{12}) into (\ref{1}) we pass to the equation
\be
\Big(\frac{d g}{d\eta}\Big)^2=F(g,q),\label{10}
 \ee
where
\begin{align}
&\begin{aligned}
&F(g,q)=3g^2-2\frac{1}{2+r}g^{2+r}-2\frac{3-q}{2-r}g^{2-r}+\frac{1-q}{1-r}g^{2-2r}-C(q),\\
&C(q)=\frac{r\big(3r^2-q(2+r)\big)}{(1-r)(4-r^2)}.\end{aligned}\label{9}
\end{align}
The right-hand side $F$ has three roots. One of them, $g=1$, corresponds to the behavior of $\om$ at infinity. The other,
$g_0<1$ and $g_1>1$, are  associated with the condition $\om(0,A)=1$. It's easy to verify that
\be
g_0\in(0,1)\; \Longleftrightarrow \; C(q)>0,\quad F''_{gg}\big|_{g=1}>0\; \Longleftrightarrow \; q>0.\label{9a}
\ee
Assuming the fulfillment of the assumptions (\ref{9a}) we obtain the function $F(g,q)$ like depicted on Fig 1.
\begin{figure}[H]\label{fig0}
\centering
\includegraphics[width=8cm]{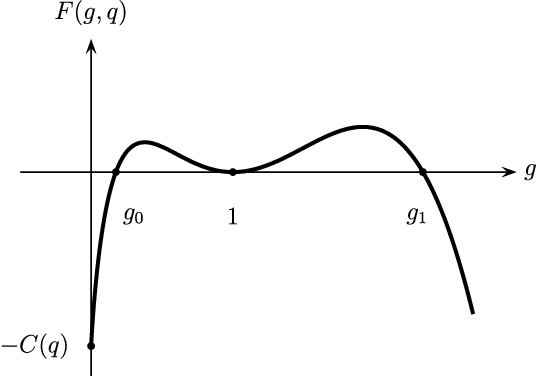}
\caption{Right-hand side of the equation (\ref{10}) in the case $r=1/2$, $q\approx0.148$. Here $g_0\approx0.175$, $g_1\approx2.455$, and $C(q)\approx1.964$.}
\end{figure}
Now we note that, for the representation (\ref{12}), the equality
\be
1-g_k^r\ipd p_k=c_3A/(\gamma+\alpha^2V),\quad k=0,1, \label{11}	
\ee
can be realized if and only if $k=0$ for $A>0$ and $k=1$ for $A<0$. In turn, equalities (\ref{11}) allow us to determine the relationship between the speed $V=V_k$ and root $g_k$. For $\alpha>0$ we have
\be
V_k=\alpha^{-2}\big\{c_3A/p_k-\gamma\big\},\quad k=0,1. \label{11a}	
\ee
Therefore, we obtain the equations  for the roots $g_k=g_k(A)$
\be
F(g_k,q)\big|_{q=q(V_k)}=0,\quad k=0,1.\label{11b}
\ee
Now we can state the Cauchy problems
\be
\frac{d g}{d\eta}=\sqrt{F(g,q)}\big|_{q=q(V_0)},\quad \eta\in(0,\infty);\quad
 g|_{\eta=0}=g_{o},\label{10a}
 \ee
 if $A>0$ and
 \be
\frac{d g}{d\eta}=-\sqrt{F(g,q)}\big|_{q=q(V_1)},\quad \eta\in(0,\infty);\quad
 g|_{\eta=0}=g_{1},\label{10b}
 \ee
 if $A<0$.

Considering in the same manner the case $\alpha=0$, we find the root $g_0=\bar{g_0}(A)$ of $F$ and the associated velocity $\bar{V}$
\be
\bar{g_0}(A)=(1-c_3A/\gamma)^{1/r}, \quad \bar{V}=c_0+c_1\gamma^2q\big(\bar{g_0}(A)\big)/c_3^2,\quad \alpha=0. \label{11c}	
\ee
Thus, taking into account the restrictions (\ref{9a}), we come to the following statement
\begin{teo}{\cite{OmNo2}}
Let the amplitude $A$ satisfy the conditions
\begin{align}
 &A\in(A_0^*,A_0^-)\bigcup(A_0^+,\infty)\quad\text{if}\quad\alpha>0,\quad\gamma_\alpha>0,\quad\text{and}\quad c_3^2>4\xi\gamma_\alpha,\label{13}\\
&A>A_0^*,\quad A\neq \overline{A}_0^\pm\quad\text{if}\quad\alpha>0,\quad\gamma_\alpha>0,\quad\text{and}\quad c_3^2=4\xi\gamma_\alpha,\label{14}\\
&A>A_0^*\quad\text{if}\quad\alpha>0,\quad\gamma_\alpha>0,\quad\text{and}\quad c_3^2<4\xi\gamma_\alpha,\label{14a}\\
&A<p_1\gamma_\alpha c_3<0\quad\text{if}\quad \alpha>0,\label{14c}
\end{align}
either $A<\gamma/c3$, $A\neq0$ if $\alpha=0$.  Then the equation (\ref{1}) for $n=3$ has the soliton solution (\ref{15}) with the velocity (\ref{11a}) if $\alpha>0$ and (\ref{11c}) if $\alpha=0$. The function $\omega(\eta,A)$ vanishes at the exponential rate,
$\omega(\eta,A)\sim \exp(-\sqrt{rq}\eta)$ for $\eta>>1$.
Here $\gamma_\alpha=\gamma+\alpha^2c_0$, $\xi=3r^2\ap^2c_1/(2+r)$, $A_0^*=p_0\gm_{\ap}/c_3$, $\overline{A}_0^\pm=p_0c_3/2\xi$, and $A_0^{\pm}=p_0\big(c_3\pm\sqrt{c_3^2-4\xi\gm_{\ap}}\big)/2\xi$.
\end{teo}
\section{Finite difference scheme}
The actual numerical simulation for the Cauchy problem for the gmKdV equation (\ref{1}) is realized for $n=3$ and for a bounded $x$-interval, $x \in [0, L]$.
For this reason we simulate the Cauchy problem by the following
mixed problem:
\begin{align}  \frac{\pa }{\pa t}\Big\{u &-\alpha^2\ve^2\frac{\pa^2 u}{\pa x^2}\Big\}+\frac{\pa}{\pa x}\Big\{c_0u+c_1u^3-c_2\ve^2\Big(\frac{\pa u}{\pa x}\Big)^2+\gamma\ve^2\frac{\pa^2 u}{\pa x^2}\notag\\
&-c_3\ve^2u\frac{\pa^2 u}{\pa x^2}\Big\}=0,  \quad x \in (0, L), \quad t \in (0, T),\label{1a}\\
&u\big|_{x=0} = u\big|_{x=L} =u_x\big|_{x=L}=0,\quad u\big|_{t=0}
= u^0 (x/\varepsilon), \label{19}
\end{align}
 where $L$, $T$, and sufficiently smooth function  $u^0$
are such that uniformly in $t \in [0,T]$
\be
\big|u(x,t)|_{x\in[0,\delta]}\big| \leq c\varepsilon^2<<1, \quad
\big|u(x,t)|_{x\in[L-\delta,L]} \big| \leq c\varepsilon^2<<1\label{20}
\ee
for a sufficiently small $\delta> 0$.

   When modeling interaction phenomena according to (\ref{1a}), (\ref{19}), we take into account that explicit formulas for multisoliton solutions to the equation (\ref{1}) remain unknown. However, solitons are functions that decay with exponential rates, so a linear combination of solitary waves (\ref{15}) approximates the exact multisoliton solution if the distances between the waves are large enough. In particular, to approximate the two-soliton solution we set
\be
u^0 =\sum_{i=1}^2A_i\om\big(\beta_i(x-x^0_i)/\varepsilon,A_i\big),\label{200}
\ee
and assume that
 \be
\frac1\ve|x^0_2-x^0_1|\geq c\varepsilon^{-\mu}>>1,\quad \mu>0.     \label{201}
\ee
Accordingly, in the case of (\ref{200}) the parameters  $L$, $T$, $\ve$ and the initial positions of the fronts $x_i^0$ should be such that the intersection points of the trajectories of solitary waves belong to the domain $Q_T = (0, L)\times(0, T)$.

Obviously, it is impossible to create any finite difference scheme
for a problem with singular perturbations  which remains reasonable
uniformly in $\varepsilon\rightarrow 0$  and $t \in (0, T)$, $T =
\mbox{const}$. So we treat $\varepsilon$ as a small but fixed
constant. However, we  fix any relation between $\varepsilon$
and finite difference scheme parameters.

To create a finite differences scheme for the equation
(\ref{1a}) we should choose appropriate approximations for the
differential and  nonlinear terms. Let us do it
separately.
\subsection{Preliminary nonlinear scheme} \
As usually, we define a mesh $Q_{T,\tau,h} = \{(x_i, t_j) \ipd (ih,
j\tau), i = 0,\dots,I$, $j = 0,\dots,J\}$ over $Q_T$ and denote
$$ y_i^j \ipd u(x_i,t_j), \quad y_{ix}^j \ipd \pa_xy_{i}^j\ipd \frac{y_{i+1}^j- y_i^j}{h},
 \quad y_{i\bar{x}}^j \ipd \pa_{\bar{x}}y_{i}^j\ipd  \frac{y_i^j - y_{i-1}^j}{h},$$
$$ y_{i\dot{x}}^j \ipd \pa_{\dot{x}}y_{i}^j\ipd \frac12(y_{ix}^j+y_{i\bar{x}}^j),
 \quad  y_{i\bar{t}}^j \ipd \pa_{\bar{t}}y_{i}^j\ipd  \frac{y_i^{j} - y_i^{j-1}}{\tau},
\quad y_{ix\bar{x}}^j = (y_{ix}^j)_{\bar{x}}.$$
Let us consider the following system of nonlinear equations:
\begin{align}
y_{i\bar{t}}^j  -\alpha^2\ve^2y_{ix\bar{x}\bar{t}}^j &+c_0y_{i\dot{x}}^j+c_1Q_1(y_{i}^j)+\varepsilon^2\big(\gamma_hy_{ix\bar{x}\dot{x}}^j +\gamma h y_{ix\bar{x}x}^{j}\big)\notag\\
& -\ve^2Q_2(y_{i}^j) = 0,\quad i = 1,\dots,I-1, \quad j = 1,2,\dots,J,\label{20a}
\end{align}
\begin{equation} \label{20b}
y_l^j = 0, \quad y_{I-l}^j = 0, \quad l=0,1,2,\quad j = 1, 2,\dots,J,
\end{equation}
\begin{equation}y_i^0 = \widetilde{u}^0 (x_i/\varepsilon), \quad i =
0,\dots,I,\label{20c}\end{equation}
where
\begin{align}
&Q_1(y)=\frac{1}{2}\{y^2y_{\dot{x}}+y(y^2)_{\dot{x}}+(y^3)_{\dot{x}}\},\quad\gamma_h=\gamma(1-h), \label{21}\\ &Q_2(y)=\pa_{\dot{x}}\Big\{c_2y_xy_{\bar{x}}+\frac{c_3}{2}\big(2yy_{x\bar{x}}+(y_x)^2-2y_xy_{\bar{x}}+(y_{\bar{x}})^2\big)\Big\},\label{22}
\end{align}
$$\widetilde{u}^0(x_l/\varepsilon)=\widetilde{u}^0(x_{I-l}/\varepsilon)=0
\quad\text{for}\quad l=0,1,2,$$
$$\widetilde{u}^0(x_i/\varepsilon)=\frac1h\int_{x_i-h/2}^{x_i+h/2}u^0
\left(\frac{\eta}{\varepsilon}\right)d\eta,\quad i=3,\dots,I-3.$$
The main properties of the terms $Q_l(y)$ are the following
\begin{align}
&h\sum_i Q_l(y_i)=0\quad l=1,2,\quad h\sum_i y_iQ_1(y_i)=0,  \label{21a}\\
&h\sum_i y_iQ_2(y_i)=\frac12(c_3-2c_2)h\sum_i y_{ix}y_{i\bar {x}}y_{i\dot{x}}.\label{22a}
\end{align}
The sense of the term $h y_{ix\bar{x}x}^{j}$ is similar to the
parabolic regularization of the gDP equation (see below). It is clear also
that  the local approximation accuracy of (\ref{20a}) is $O(\tau +
h^2)$ for sufficiently smooth solution. It should be noted also that  a similar approach to the nonlinearity $Q_1$ digitization has been presented and successfully used in \cite{Sep, PSV}.

To simplify the notation, we  write
$$y \ipd y_i^j,  \quad \check{y} \ipd y_i^{j-1}.$$
So, the short form of the equation (\ref{20a}) is the following:
\begin{equation} \label{23}
y_{\bar{t}} -\alpha^2\ve^2y_{x\bar{x}\bar{t}}+c_0y_{\dot{x}}+c_1Q_1(y)+\varepsilon^2\big(\gamma_hy_{x\bar{x}\dot{x}} +\gamma h y_{x\bar{x}x}\big)
 -\ve^2Q_2(y) = 0.
\end{equation}

Our first result consists of  obtaining of discrete analogs of the
equalities (\ref{4}) and (\ref{3}). Multiplying (\ref{20a}) by $1$ and
$y$, summing over $i$, and using some trivial equalities (see \cite{OmNo1}) it is easy to
obtain the following
\begin{lem} \label{lem31}
Let the nonlinear system of algebraic equations
(\ref{20a})-(\ref{20c}) have the unique solution $y_{i}^j$, $i =
1,...,I-1$, $j = 1,2,...,J$, and let $\ve>0$ be a constant. Then
uniformly in $j \leq J$ the following relations hold:
\begin{align}
 &\pa_{\bar{t}}\,h\sum y^j=0,\label{24}\\
&\pa_{\bar{t}}\Big\{\|y^j\|^2+\alpha^2\|\ve y_x^j\|^2\Big\}+\tau\Big\{\|y_{\bar{t}}^j\|^2+\alpha^2\|\ve y_{x\bar{t}}^j\|^2\Big\}+\gamma h^2\|\ve y_{x\bar{x}} ^j\|^2\notag\\
&\qquad\qquad\qquad=\ve^2(c_3-2c_2)h\sum y_{x} ^jy_{\bar{x}}^jy_{\dot{x}} ^j,\label{25}
\end{align}
\end{lem}
Here and in what follows $\sum$ denotes the summation over all $i$ and $\| \cdot \|$ is the discrete version of
the  ${L}^2(0,L)$ norm, namely
\be
\sum f=\sum_{i=1}^{I-1}f_i\,,\quad\| f \|^2 =h\sum_{i=1}^{I-1}|f_i|^2.\label{28}
\ee
Next, we note that the equality (\ref{25}), as well as (\ref{4}), does't imply any regularity of the solution if $c_3\neq2c_2$. However, there exists a special case, when (\ref{25}) allows us to analyze the equation (\ref{23}) solution.
\begin{lem} \label{lem32} Let $T>0$ be a  constant, $\alpha>0$, and $u^0,\;u_x^0\in L^2(0,L)$.
Moreover, let for each $j= 0,1,2,...,J$ there exist $k(j)$ such that
\begin{equation} \label{29}
y^j_{k(j)+i}=y^j_{k(j)-i},\quad i=0,1,\dots,
\end{equation}
where $y^j_{i}=0$ for $i\leq0$ and $i\geq I$. Then under the assumptions of Lemma \ref{lem31}
\begin{align}
\| y^j\|^2 +\alpha^2\|\ve y_x^j\|^2&+\tau \Big\{\|| y_{\bar{t}} \||^2(j)+\alpha^2\||\ve y_{x\bar{t}}^j\||^2(j)\Big\}\notag\\
& +\gamma h^2\||\ve  y_{x\bar{x}}\||^2(j) =\| y^0\|^2+\alpha^2\|\ve y_{x}^0\|^2.\label{30}
\end{align}
\end{lem}
Here and in what follows $\|| \cdot \||(j)$
is the discrete version of the
${L}^2\big((0,L)\times(0,t_j)\big)$ norm, namely
\be\|| f \||^2(j)
= \tau \sum_{k=1}^j \| f^k \|^2.\label{31}\ee
As a consequence of this lemma one can prove a convergence result.
Namely, let $y_{\tau,h}(x,t)$ be an extension of the net-function
$y_i^j$ which satisfies the same estimate (\ref{30}) and the evenness assumption in the sense of (\ref{29}) (see e.g.
\cite{Lad}), and let $u(x,t)$ denote the solution of the problem
(\ref{1a}), (\ref{19}). Then, using the standard technic (see
e.g. \cite{Li}), one can prove the theorem
\begin{teo} \label{teo31}
Let the assumptions of Lemma \ref{lem32} be satisfied. Then there
exists a subsequence $y_{\bar{\tau},\bar{h}}(x,t)$ such that
\be\label{32a}
y_{\bar{\tau},\,\bar{h}} \to u\quad\text{$\ast$-weakly in}\quad
L^\infty\big((0,T);W_2^1(0,L)\big)\big)
\ee
as $\tau$, $h\to0$, where $W_2^l$ denotes the Sobolev space.
\end{teo}
 Similar to (\ref{30}) one can obtain stronger estimates.
\begin{lem} \label{lem33}
Under the assumption of Lemma \ref{lem32} let  $u^0\in W_2^2$. Then uniformly in $j$:
\begin{align}
&\|\ve y_x^j\|^2 +\alpha^2\|\ve^2 y_{x\bar{x}}^j\|^2+\tau \Big\{\||\ve y_{x\bar{t}} \||^2(j)+\alpha^2\||\ve^2 y_{x\bar{x}\bar{t}}\||^2(j)\Big\}\notag\\
& \qquad\qquad+\gamma h^2 \||\ve^2 y_{x\bar{x}\dot{x}}\||^2(j) =\|\ve y_x^0\|^2+\alpha^2\|\ve^2 y_{x\bar{x}}^0\|^2,\label{32}\\
&\|y_{\bar{t}}^j\|^2+\|\ve y_{x\bar{t}}^j\|^2
 \leq \ve^{-4}C\big(\|y^0\|,\|\ve y_{x}^0\|,\|\ve^2 y_{x\bar{x}}^0\|\big),\label{33}
\end{align}
where $C(v,w,z)$ does not depend on
$\tau$,  $h$, and $\ve$.
\end{lem}
\subsection{Linearization}\
Now we should verify the solvability of the equation (\ref{23})
for any fixed $j \geq 1$, that is of the system of nonlinear equations
\begin{align}
y -\alpha^2\ve^2y_{x\bar{x}}+\tau\Big\{c_0y_{\dot{x}}+c_1Q_1(y)
&+\varepsilon^2\big(\gamma_hy_{x\bar{x}\dot{x}} +\gamma h y_{x\bar{x}x}\big)\notag\\
& -\ve^2Q_2(y)\Big\}= \check{y}-\alpha^2\ve^2\check{y}_{x\bar{x}},\label{37}
\end{align}
as well as select a way to linearize the nonlinear terms. To this aim
let us construct the sequence of vector functions $\varphi(s) \ipd
\{\varphi_0(s), . . ., \varphi_I(s)\}$, $s \geq 0$, such that
$\varphi(0) = \check{y}$ and $\varphi\ipd\varphi(s)$ for $s \geq 1$ satisfies
the equation
\begin{align}
&\varphi -\alpha^2\ve^2\varphi_{x\bar{x}} +\tau\big\{c_0\varphi_{\dot{x}}+c_1\textit{R}_1(\bar{\varphi},\varphi)+\varepsilon^2\big(\gamma_h\varphi_{x\bar{x}\dot{x}} +\gamma h \varphi_{x\bar{x}x}\big)\notag\\
&-\ve^2\textit{R}_2(\bar{\varphi},\varphi)\big\}= \check{y}-\alpha^2\ve^2\check{y}_{x\bar{x}}
+\tau\big\{2c_1Q_1(\bar{\varphi})-\ve^2Q_2(\bar{\varphi})\big\}, \label{38}\\
&\varphi_l = 0, \quad\varphi_{I-l} = 0, \quad l=0,1,2,\notag
\end{align}
where  $\bar{\varphi} \ipd \varphi(s-1)$. To obtain (\ref{38}) we use the identity
\begin{align}
&Q_1(\varphi)=Q_1(\bar{\varphi}+\mathrm{w})=Q_1(\bar{\varphi})+\textit{R}_1(\bar{\varphi},\mathrm{w})
+\textit{R}_1(\mathrm{w},\bar{\varphi})+Q_1(\wm),\label{39}\\
&Q_2(\varphi)=Q_2(\bar{\varphi}+\mathrm{w})=Q_2(\bar{\varphi})+\textit{R}_2(\bar{\varphi},\mathrm{w})+Q_2(\wm),\label{39a}\\
&\textit{R}_1(u,v)=\frac{1}{2}\big\{u^2v_{\dot{x}}+2u(uv)_{\dot{x}}+3(u^2v)_{\dot{x}}+2uvu_{\dot{x}}+v(u^2)_{\dot{x}}\big\},\label{40}\\
&\textit{R}_2(\bar{\varphi},\mathrm{w})=\pa_{\dot{x}}\Big\{(c_2-c_3)\big(\bar{\varphi}_{x}\mathrm{w}_{\bar{x}}+\bar{\varphi}_{\bar{x}}\mathrm{w}_x\big)\notag\\
&\qquad\qquad+c_3\big(\bar{\varphi}\mathrm{w}_{x\bar{x}}+\bar{\varphi}_{x}\mathrm{w}_{x}+\bar{\varphi}_{\bar{x}}\mathrm{w}_{\bar{x}}+\bar{\varphi}_{x\bar{x}}\mathrm{w}\big)\Big\},\label{41}
\end{align}
where $\wm=\vp-\bar{\varphi}$. Next we neglect the quadratic in $\mathrm{w}$ terms in (\ref{39}), (\ref{39a}), and note that
\be \label{42}
\textit{R}_1(\bar{\varphi},\mathrm{w})=\textit{R}_1(\bar{\varphi},\varphi)-3Q_1(\bar{\varphi}),\quad
\textit{R}_2(\bar{\varphi},\mathrm{w})=\textit{R}_2(\bar{\varphi},\varphi)-2Q_2(\bar{\varphi}).
\ee
 Continuing, we  we note that the solvability of the algebraic system
(\ref{38}) is obvious for sufficiently small $\tau/h^3$. In order to estimate $\|\vp\|$ let us use the identities (\ref{39}) - (\ref{42}) again and rewrite (\ref{38}) as
follows:
\begin{align}
&\varphi -\alpha^2\ve^2\varphi_{x\bar{x}} +\tau\big\{c_0\varphi_{\dot{x}}+c_1Q_1(\varphi)+\varepsilon^2\big(\gamma_h\varphi_{x\bar{x}\dot{x}} +\gamma h \varphi_{x\bar{x}x}\big)\notag\\
&-\ve^2Q_2(\varphi)\big\}= \check{y}-\alpha^2\ve^2\check{y}_{x\bar{x}}
+\tau\big\{c_1\textit{R}_1(\mathrm{w},\bar{\varphi})+c_1Q_1(\mathrm{w})-\ve^2Q_2(\mathrm{w})\big\}. \label{43}
\end{align}
Similarly to (\ref{25}) the equation (\ref{43}) yields
\begin{align}
&\|\varphi\|^2+\alpha^2\|\ve \varphi_x\|^2+\tau\gamma h^2\| \ve\varphi_{x\bar{x}}\|^2-\ve^2(c_3/2-c_2)h\sum \varphi_{x}\varphi_{\bar{x}}\varphi_{\dot{x}}\notag\\
&=h\sum \varphi\{\check{y}-\alpha^2\ve^2\check{y}_{x\bar{x}}\}+\tau h\sum \varphi\{c_1\textit{R}_1(\mathrm{w},\bar{\varphi})+c_1Q_1(\mathrm{w})-\ve^2Q_2(\mathrm{w})\}.\label{44}
\end{align}
Again we should assume the existence of a special solution $\varphi$, which satisfies the evenness condition (\ref{29}).
The next step is the estimation of the  discrepancy $\mathrm{w}$.
Assuming the existence of the special even solution $y^k$ of
the equation  (\ref{23}) with $k = 1, 2,..., j-1$, and subtracting one equation (\ref{38}) (for $s=s_0-1$) from the
another one (for $s=s_0$), we obtain:
\begin{align}
&\mathrm{w}-\alpha^2\ve^2\mathrm{w}_{x\bar{x}} +\tau \big\{c_0\mathrm{w}_{\dot{x}}+c_1\textit{R}_1(\check{y},\mathrm{w})+\varepsilon^2\big(\gamma_h\mathrm{w}_{x\bar{x}\dot{x}} +\gamma h \mathrm{w}_{x\bar{x}x}\big)\notag\\
&-\ve^2\textit{R}_2(\check{y},\mathrm{w})\big\}=\tau(\check{y}_{\bar{t}}-\alpha^2\ve^2\check{y}_{x\bar{x}\bar{t}})\quad\text{for}\quad s=1,\label{46}\\
&\mathrm{w}-\alpha^2\ve^2\mathrm{w}_{x\bar{x}} +\tau \big\{c_0\mathrm{w}_{\dot{x}}+c_1\textit{R}_1(\bar{\varphi},\mathrm{w})
+c_1\textit{R}_1(\bar{\mathrm{w}},\bar{\bar{\varphi}})+\varepsilon^2\big(\gamma_h\mathrm{w}_{x\bar{x}\dot{x}} +\gamma h \mathrm{w}_{x\bar{x}x}\big)\notag\\
&-\ve^2\textit{R}_2(\bar{\varphi},\mathrm{w})\big\}=-\tau\big(c_1Q_1(\bar{\mathrm{w}})-\ve^2Q_2(\bar{\mathrm{w}})\big)\quad\text{for}\quad s>1, \label{47}
\end{align}
where $\bar{\mathrm{w}}\ipd\varphi(s-1)-\varphi(s-2)$.
Applying the standard techniques we verify the following estimates for
$\varphi$ and $\mathrm{w}$ (for the proof see Attachment):
\begin{lem} \label{lem36}
Let the assumptions of Lemma \ref{lem33} be satisfied and let $\varphi$ satisfy the evenness condition (\ref{29}).
Suppose also that
\be \tau\leq q_1\ve h^2,\quad h\leq q_2\ve,\label{48} \ee
where constants $q_i>0$ are sufficiently small.
Then
\begin{align}
&\| \varphi \|_{(2,\tau)}^2\big\{1-(\ve q_1^3q_2^4)^{1/4}\big(\|\mathrm{w}\|_{(2,\tau)}^2 +(q_1q_2)^{1/4}h^{3/4}\|\mathrm{w}\|_{(1)}^2 \big)\big\} \notag\\
&\leq c_1\|\check{y}\|^2_{(2)}+c_2 \big\{\|\mathrm{w}\|^4_{(2,\tau)}+\|\mathrm{w}\|^6_{(1)}\big\},  \label{49}\\
&\|\mathrm{w}\|_{(2,\tau)}^2 \leq c_3 \tau^{2}\|\check{y}_{\bar{t}}\|^2_{(1)}\quad\text{for}\quad s=1,\label{50}\\
&\|\mathrm{w}\|_{(2,\tau)}^2\leq c_4(\tau\ve^{-2})^2\|\bar{\mathrm{w}}\|_{(1)}^6+c_5\|\bar{\mathrm{w}}\|_{(2,\tau)}^4,\quad \text{for} \quad s> 1,\label{51}
\end{align}
where
\begin{align}
&\|f\|^2_{(1)}\ipd\|f\|^2+\|\ve f_x\|^2,\quad \|f\|^2_{(2)}\ipd\|f\|^2+\|\ve f_x\|^2+\|\ve^2f_{x\bar{x}}\|^2,\notag\\
&\quad\|f\|^2_{(2,\tau)}\ipd\|f\|^2+\|\ve f_x\|^2+\tau h^2\|\ve f_{x\bar{x}}\|^2,\label{52}
\end{align}
 $c_i > 0$ denote constants which do not depend on $h$, $\tau$, $\ve$, and $s$.
\end{lem}
In particular, estimates (\ref{48})-(\ref{51}) imply for $s=2$
\begin{align}
 &\|\mathrm{w}(2)\|^2_{(2,\tau)} \leq c\big(\tau \ve^{-2}\big)^{2}\tau^6+(c'\tau^2)^2\leq c\tau^4,\notag\\
 &\| \varphi(2) \|_{(2,\tau)}^2(1-c\tau^4)\leq c\|\check{y}\|^2_{(2)}+c'\tau^8. \label{53}
\end{align}
Collecting the estimates (\ref{33}) and (\ref{48})-(\ref{51}) together, we finally reach the conclusion that
the terms of the $\mathrm{w}$-sequence vanish very rapidly,
\begin{equation} \label{53a}
 \|\mathrm{w}(s)\|^2_{(2,\tau)} \leq \big(c'\tau\big)^{2^s}.
\end{equation}
 Note next that by virtue of (\ref{49})-(\ref{51}), the terms of
$\varphi$-sequence are bounded uniformly in $s$
\begin{equation} \label{55}
\|\varphi(s)\|^2_{(2,\ve)} \leq c\|y^j\|^2+O(\tau^4).
\end{equation}
Furthermore, for any $n > 0$
\begin{align}\|\varphi(s+n) -
\varphi(s)\|_{(2,\ve)} &\leq
\sum_{i=1}^n\|\mathrm{w}_{s+i}\|_{(2,\ve)} \notag\\
&\leq \|\mathrm{w}_{s+1}\|_{(2,\ve)} \sum_{i=1}^\infty
\frac{\|\mathrm{w}_{s+i}\|_{(2,\ve)}}{\|\mathrm{w}_{s+1}\|_{(2,\ve)}}
\leq c \|\mathrm{w}_{s+1}\|_{(2,\ve)}.\notag
\end{align}
This implies the main statement of this subsection:
\begin{teo} \label{teo33}
Let the assumption Lemma \ref{lem36} be satisfied. Then  the
sequence $\varphi(s)$ converges in the $H^2_{\tau,h,\ve}$ sense to the
solution of the equation (\ref{37}). Moreover,
\begin{equation} \label{56}
\|y - \varphi(2)\| \leq c\tau^4,
\end{equation}
where $H^2_{\tau,h,\ve}$ is the space with the norm (\ref{52}) and a constant $c > 0$  dos not depend on $h$, $\tau$, and $\ve$.
\end{teo}
\subsection{Numerical simulation}\
To solve the system of linear equations  (\ref{38})
we use the Gauss method adapted to systems with five non-zero
diagonals. This implies the efficiency of the scheme in the sense
that it executes $O(I)$ arithmetic operations to pass to the next
time-level. In accordance with (\ref{56}) we stop the iterative calculation of $\varphi(s)$ at the second step setting  $y^{j}=\varphi(2)$. Clearly, this implies the appearance of an error. However,  our results of numerical simulations justify this decision (see below).

In order to define soliton initial data  we solve numerically the equation (\ref{11b}) and define $g(0)=g_*$, where $g_*=g_0$ for $A>0$ and $g_*=g_1$ for $A<0$.  Next, to avoid the non uniqueness in the problems (\ref{10a}), (\ref{10b}) we calculate
$$
g^*(h)=g_*+\frac{1}{4}h^2\frac{dF}{dg}\big|_{g=g_*}+\frac{1}{4}\frac{h^4}{4!}\frac{dF}{dg}\frac{d^2F}{dg^2}\big|_{g=g_*}+\left.\frac{1}{4}\frac{h^6}{6!}\frac{dF}{dg}\corc{3\frac{d^3F}{dg^3}+\pare{\frac{d^2F}{dg^2}}^2}\right|_{g=g_*}
$$
and solve the similar (\ref{10a}) (if $A>0$ or (\ref{10b}) if $A<0$) problem
\be
\frac{d g}{d\eta}=\sqrt{F(g,q)},\quad \eta\in(h,\infty);\quad
 g|_{\eta=h}=g^*(h),\label{57}
 \ee
using the fourth order Runge-Kutta method. The last step is the determination of the profile $\om(\eta,A)$ in accordance with the rule (\ref{12}).

\textbf{Example 1.}
 When $\alpha=c_0=c_2=c_3=0$, $c_1=2$, $\gamma=1$, and $n=3$,  (\ref{1}) is the modified KdV equation. For definiteness, here and in what follows we set $\ve=0.1$. We test the finite difference scheme (\ref{38}) by calculating the difference $Er$ between the numerical and exact mKdV solitons, and also check the fulfilment of the conservation
laws \eqref{24} and \eqref{25}, see Fig. 2 and Table 1. The motion of the numerical mKdV soliton  is shown in Fig. 3.
\begin{figure}[H]
\centering
\includegraphics[width=8cm]{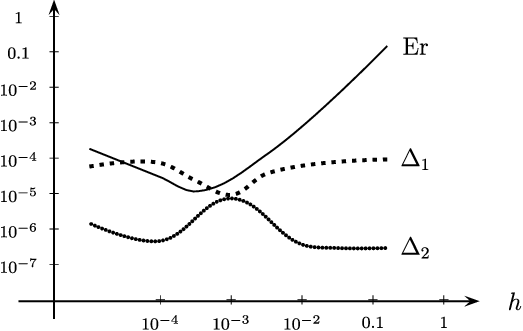}
\caption{Behavior of the error-functions $Er$ and $\Delta_k$ (\ref{58}) for the mKdV soliton with $A=1.2$ and $\tau=h^2$}
\label{f1}
\end{figure}
Here
\be
Er=\max_{i=0,\dots,I}|u_{exact}(x_i,t_j)-y_i^j|,\quad
\Delta_k=\max_{j'=0,\dots,j}|E_{k}^{j'}-E_{k}^{0}|,\quad k=1,2,\label{58}
\ee
$u_{\text{exact}}=A\cosh^{-1}(\beta(x-Vt-x^0)/\ve)$ is the exact mKdV soliton, $V=A^2$, $\beta=A$, and $y_i^j$ is the  numerical wave  at $x_i=x^{0}+ih$, $t_j=j\tau$.
 The energies $E_{k}^{j}$ are calculated in accordance with formulas (\ref{21a}), (\ref{22a}),
\be
\begin{aligned}
&E_{1}^{j}=h\sum_{i=0}^{I} y_i^j,\quad E_{2}^{j}=\|y^j\|^2+\alpha^2\|\ve y_x^j\|^2+\tau^2\sum_{j'=1}^{j}\Big\{\|y_{\bar{t}}^{j'}\|^2+\alpha^2\|\ve y_{x\bar{t}}^{j'}\|^2\Big\}\notag\\
&\hspace*{2cm}+\gamma h^2\tau\sum_{j'=1}^{j}\|\ve y_{x\bar{x}}^j\|^2+\ve^2(2c_2-c_3)\tau h\sum_{j'=1}^{j}\sum_{i=0}^{I} y_{ix}^jy_{i\bar{x}}^jy_{i\dot{x}}^j.\notag
\end{aligned}\ee
\begin{table}[H]
\begin{center}
\begin{tabular}{l|l|l|l|l|l|l|l|l|l|l}
$h\pare{\times10^{-3}}$&$20.0$&$12.5$&$10.0$&$7.1$&$5.5$&$5.0$&$4.5$&$4.1$&$3.8$&$3.3$\\ \hline
$\E\pare{\times10^{-3}}$&$170.7$&$71.4$&$46.5$&$24.0$&$14.6$&$1.8$&$0.9$&$0.8$&$7.0$&$9.0$\\ \hline
$\D_1\pare{\times10^{-6}}$&$18.3$&$14.4$&$11.1$&$6.5$&$3.9$&$4.1$&$4.3$&$2.9$&$4.6$&$9.0$\\ \hline
$\D_2\pare{\times10^{-6}}$&$0.5$&$0.4$&$0.5$&$0.6$&$0.6$&$1.9$&$1.1$&$3.5$&$0.7$&$1.2$\\ \hline
\end{tabular}
\caption{The errors $\E$, $\D_k$ for $h\in [3.3\times 10^{-3}, 0.2]$ and $A=1.2$ at time $t=t_j=1$.}
\label{t1}
\end{center}\end{table}
\begin{figure}[H]
\centering
\includegraphics[width=14cm]{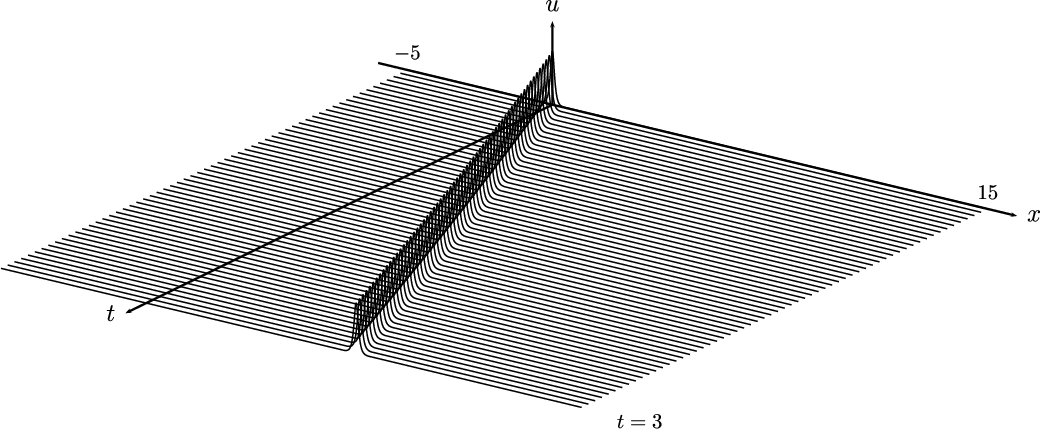}
\caption{Dynamics of the mKdV soliton with $A=1.6$, $h=4.1\times 10^{-3}$, and $\tau=h^2$}
\end{figure}
\textbf{Example 2.}
Let us consider mgDP equation (\ref{1}) in the case
 \be
 c_3=c_2=2,\quad \gamma=2,\quad c_1= c_0=\alpha=1.\label{60a}
\ee
Then  $r=1/2$ and  Theorem 1 guarantees the existence of solitons with amplitudes under the assumption \eqref{13}, where $A_0^*=0.33$, $A_0^-=1.9$, and $A_0^+=2.55$.
\begin{table}[H]
\begin{center}
\begin{tabular}{l|l|l|l|l|l|l|l|l|l|l}
\scriptsize{$h\pare{\times10^{-3}}$}&\scriptsize{$20.0$}&\scriptsize{$12.5$}&\scriptsize{$10.0$}&\scriptsize{$7.1$}&\scriptsize{$5.5$}&\scriptsize{$5.0$}&\scriptsize{$4.5$}&\scriptsize{$4.1$}&\scriptsize{$3.8$}&\scriptsize{$3.3$}\\ \hline
\scriptsize{$\D_1\pare{\times10^{-5}}$}&\scriptsize{$0.5$}&\scriptsize{$2.5$}&\scriptsize{$7.3$}&\scriptsize{$5.2$}&\scriptsize{$46.0$}&\scriptsize{$20.8$}&\scriptsize{$102.3$}&\scriptsize{$201.3$}&\scriptsize{$523.2$}&\scriptsize{$1037.3$}\\ \hline
\scriptsize{$\D_2\pare{\times10^{-3}}$}&\scriptsize{$328.2$}&\scriptsize{$246.1$}&\scriptsize{$146.1$}&\scriptsize{$97.2$}&\scriptsize{$84.1$}&\scriptsize{$38.3$}&\scriptsize{$1.1$}&\scriptsize{$3.5$}&\scriptsize{$41.5$}&\scriptsize{$317.9$}\\ \hline
\end{tabular}
\caption{The errors $\D_k$ for $h\in [3.8\times 10^{-3}, 0.2]$ and $A=1.2$ at time $t=t_j=1$.}
\label{t12}
\end{center}\end{table}
\begin{figure}[H]
\centering
\includegraphics[width=8cm]{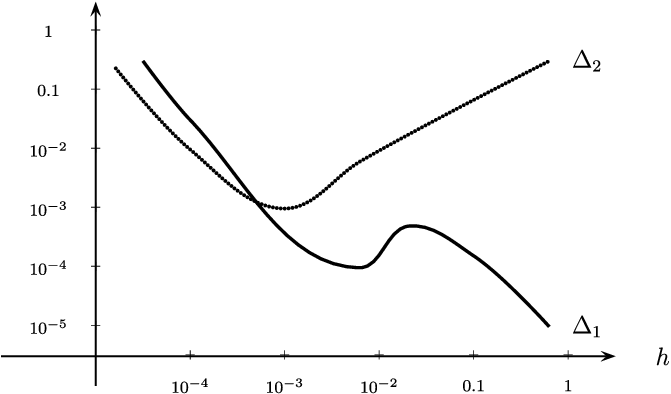}
\caption{Behavior of the error-functions $\Delta_k$ for the mgDP soliton with $A=1.2$  in the case (\ref{60a})}
\label{f2}
\end{figure}
Next, for $r=1/2$ the function $F(g,q)$ (\ref{9}) is the 5-degree polynomial
$$
F(g,q)=(1/15)(z-1)^2\big(-12z^3+21z^{2}+(20q-6)z+10q-3\big)\big|_{z=g^{1/2}}.
$$
Thus, the real root $g=g^*$ can be found analytically using the Cardano's formula. We set $A=1.2$ and define the initial data
by solving the problem (\ref{57}). The graph in Fig. 4 and Table 2 confirm the stability of the wave propagation.
Figures 5 and 6 depict the evolution of one and two solitons respectively for the case (\ref{60a}).
\begin{figure}[H]
\centering
\includegraphics[width=12cm]{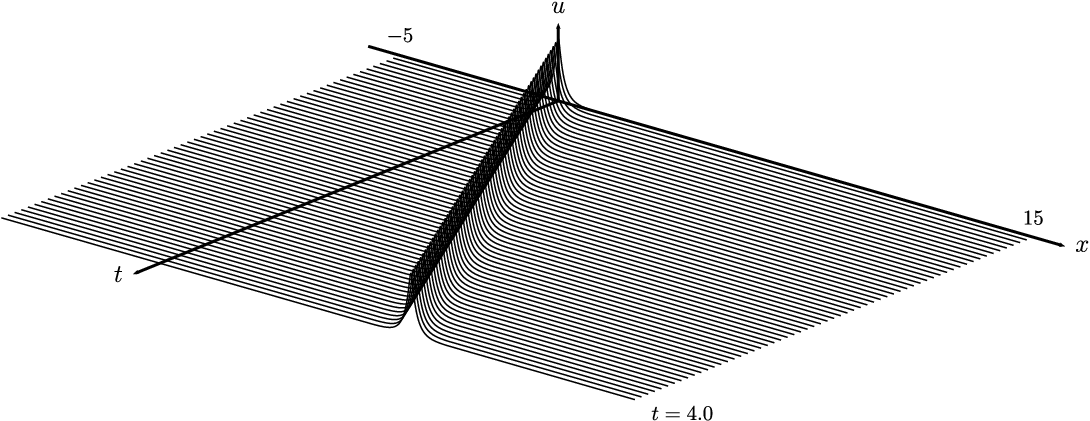}
\caption{Evolution of the mgDP soliton with $A=1.2$ in the case (\ref{60a})}
\label{f3}
\end{figure}

\begin{figure}[H]
\centering
\includegraphics[width=12cm]{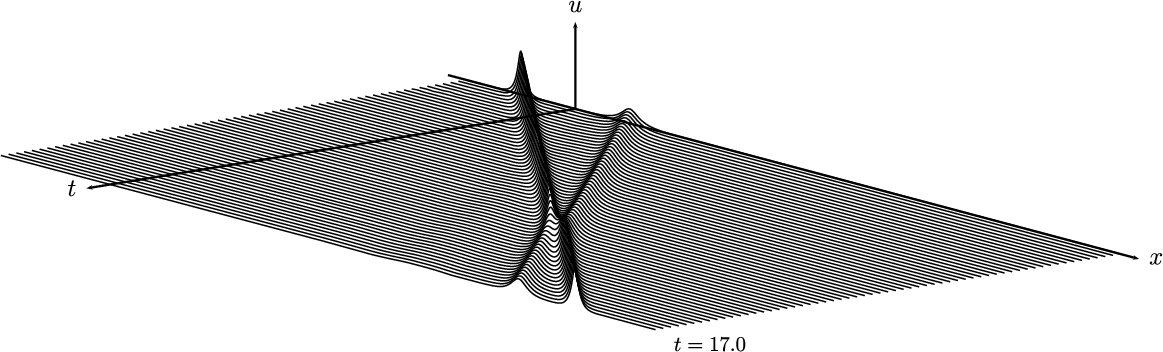}
\caption{Collision of two solitons with $A_1=1.2$, $A_2=0.5$ in the case (\ref{60a})}
\label{f4}
\end{figure}
\textbf{Example 3.}
Let us now consider mgDP equation (\ref{1}) in the case
 \be
 \ds \ap=c_1=c_3=c_2=1,\quad \gamma=c_0=2.\label{61a}
\ee
Then  $r=1/2$ and  Theorem 1 guarantees the existence of solitons with amplitudes under the assumption \eqref{14a}, where $A_0^*=0.20$.
\begin{table}[H]
\begin{center}
\begin{tabular}{l|l|l|l|l|l|l|l|l|l}
\scriptsize{$h\pare{\times10^{-3}}$}&$6.2$&$6.0$&$5.8$&$5.7$&$5.5$&$5.4$&$5.2$&$5.1$&$5.0$\\ \hline
\scriptsize{$\D_1\pare{\times10^{-4}}$}&$11.4$&$6.1$&$16.8$&$42.1$&$10.5$&$18.7$&$5.0$&$20.9$&$51.0$\\ \hline
\scriptsize{$\D_2\pare{\times10^{-3}}$}&$17.0$&$28.4$&$23.0$&$18.6$&$8.0$&$7.3$&$27.3$&$18.5$&$10.1$\\ \hline
\end{tabular}
\caption{The errors $\D_k$ for $h\in [5\times 10^{-3}, 6\times 10^{-3}]$ and $A=1.5$ at time $t=t_j=1$.}
\label{t13}
\end{center}\end{table}
\begin{figure}[H]
\centering
\includegraphics[width=10cm]{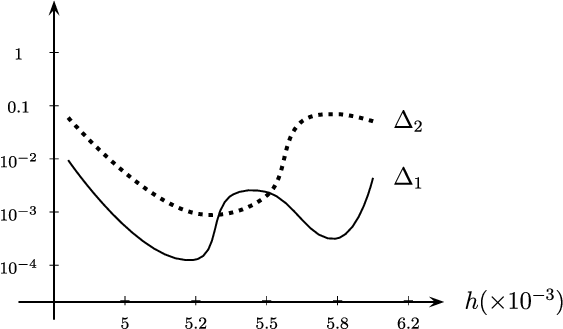}
\caption{Behavior of the error-functions $\Delta_k$ for the mgDP soliton with $A=1.5$  in the case (\ref{61a})}
\label{f23}
\end{figure}
Next, similar to Example $2$, for $r=1/2$ the function $F(g,q)$  is the 5-degree polynomial (\ref{9}). Thus, the real roots $g=g_0\in (0,1)$ and $g=g_1>1$ can be found analytically using the Cardano's formula. We use either $A>0$ or $A<0$ with $g_0$ or $g_1$ respectively and define the initial data by solving the problem (\ref{57}). The graph in Fig. 8 and Table 3 confirm the stability of the antisoliton propagation. Figure 9 depicts  soliton-antisoliton interaction for the case (\ref{61a}).
\begin{figure}[H]
\centering
\includegraphics[width=12cm]{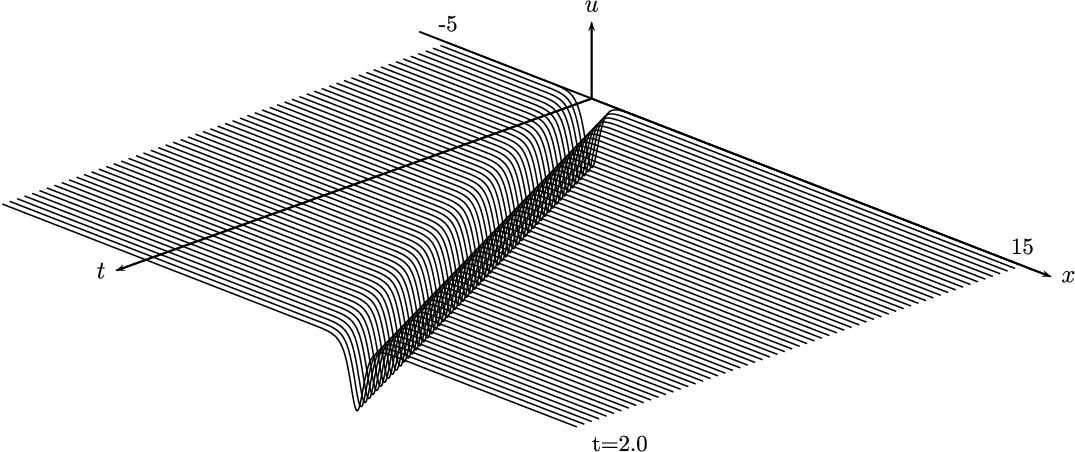}
\caption{Evolution of the mgDP antisoliton with $A=-1.8$ in the case (\ref{61a}) with $h=5.2\times 10^{-3}$ and $g_1= 1.73473808$}
\label{f33}
\end{figure}

\begin{figure}[H]
\centering
\includegraphics[width=13cm]{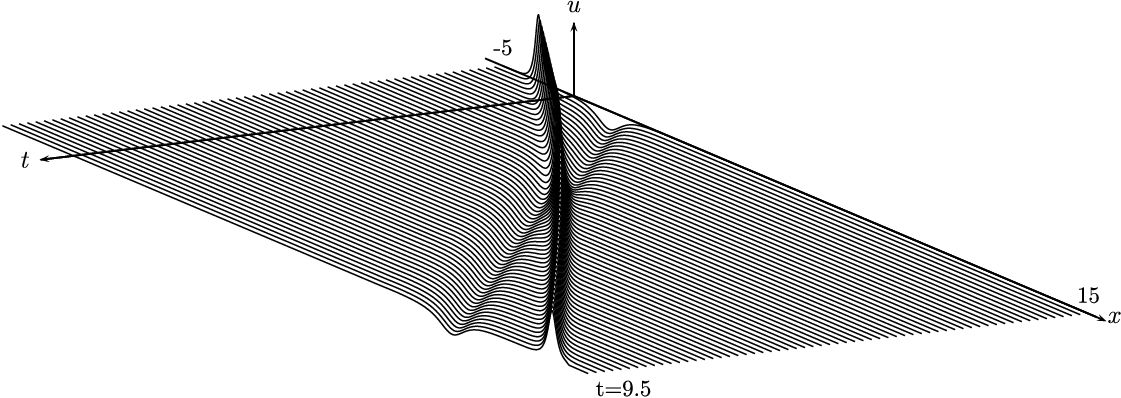}
\caption{Collision of the soliton-antisoliton with $A_1=1.8$, $A_2=-0.5$ in the case (\ref{61a})}
\label{f43}
\end{figure}
\section{Conclusion}
Energy estimates and results of numerical experiments confirm the adaptation of the balance laws (\ref{4}), (\ref{3}) for the gmKdV equation (\ref{1}) by difference scheme (\ref{38}). This implies the stability of the motion of even-shaped waves in both the analytical and numerical sense. Moreover, the scheme (ref{38}) remains stable even in the case of soliton interaction, when, generally speaking, the solution ceases to be even at the time instant of collision of waves. At the same time, it turned out that equation (\ref{1}) with $n=3$ is much more sensitive to the accuracy of  the initial data than the gDP equation with $n=2$. Indeed, the gmKdV equation requires  accuracy $O(h^8)$ in the initial data approximation for (\ref{57}), while for the gDP equation it was sufficient to use only three terms of the Taylor expansion, see \cite{OmNo1}.

{\textbf{Data availability}}

No data was used for the research described in the article.

\section{Attachment}
In what follows we use the notation
$$
||f||_p=\Big(h\sum_{i=1}^{I-1}|f_i|^p\Big)^{\frac1p}
$$
for the discrete analogs of the ${\textsf{L}}^p(0,L)$  norm.
Again, for simplicity we write
$||f||\ipd||f||_p$ if $p=2$.

Our main tools are the discrete versions of the H\"{o}lder and the
Gagliardo - Nirenberg inequalities, namely,
$$ h\left|\sum_{i=0}^N f_i g_i \right| \leq \|f\|_p \|g\|_q, \quad \frac{1}{p} +
\frac{1}{q} = 1, \quad 1 < p,q < \infty,$$
\begin{equation} \label{A1}
\| \partial_x^r f \|_p \leq c \| f \|^{1-\theta}\big\{\| f
\|^2+\|\partial_x^{\ell}f\|^2\big\}^{\theta/2}, \quad \theta \ell = \frac{1}{2} + r -
\frac{1}{p},\quad 0<\theta<1,
\end{equation}
where $c$ is a constant which does not depend on $h$. Let us
recall that the Gagliardo-Nirenberg inequality  is the
multiplicative form of the embedding theorem (see e.g. \cite{Lad}).
The proof of the discrete version of the Gagliardo-Nirenberg
inequality can be found for example in \cite{Lad}. For $x\in\R^1$
the proof of (\ref{A1}) is trivial. Recall also that for $f$ from the Sobolev space $H_0^l(0,L)$ of functions with zero
value on the boundary,
\begin{equation} \label{A2}
\max_i|f_i|\leq \sqrt2 \| f \|^{1/2} \| f_x\|^{1/2}.
\end{equation}
Furthermore, we use the well-known identities $$
(yg)_x=y_xg+yg_x+hy_xg_x,\quad
(yg)_{\bar{x}}=y_{\bar{x}}g+yg_{\bar{x}}-hy_{\bar{x}}g_{\bar{x}},$$
so that \be
(yg)_{\dot{x}}=y_{\dot{x}}g+yg_{\dot{x}}+\frac{h^2}{2}(y_xg_x)_{\bar{x}}.\label{A3}\ee
\subsection{Proof of Lemma \ref{lem33}}
Multiplying the equation
(\ref{23}) by $\ve^2y_{x\bar{x}}$, summing over $i$ and ``integrating" by parts  we obtain
\begin{align}
&\pa_{\bar{t}}\big\{\|\ve y_x\|^2 +\alpha^2\|\ve^2 y_{x\bar{x}}\|^2\big\}+\tau \big\{\|\ve y_{x\bar{x}}\|^2 +\alpha^2\|\ve^2 y_{x\bar{x}\bar{t}}\|^2\big\}\notag\\
&\qquad\qquad\;+\gamma h^2 \|\ve^2 y_{x\bar{x}\dot{x}}\|^2=c_1\ve^2h\sum Q_1y_{x\bar{x}}+\ve^4h\sum Q_2y_{x\bar{x}}.\notag
\end{align}
However, for even $y$
$$
h\sum_{i=1}^{I-1}Q_l(y_i)y_{x\bar{x}}=0,\quad l=1,2.
$$
Thus, uniformly in $j$
\begin{equation}
\|\ve y_x^j\|^2 +\alpha^2\|\ve^2 y_{x\bar{x}}^j\|^2\leq c. \label{A4}
\end{equation}
Next we multiply (\ref{23}) by $y_{\bar{t}}$. Summation over $i$ implies the inequality
$$
\| y_{\bar{t}}\|^2 +\alpha^2\|\ve y_{x\bar{t}}\|^2\leq c_0\| y_{\bar{t}}\|\|y_{x}\|+\ve^2\gamma \| y_{x\bar{t}}\|\|y_{x\bar{x}}\|
+h|\sum \{c_1Q_1+\ve^2Q_2\}y_{\bar{t}}|.
$$
Furthermore, applying the Gagliardo - Nirenberg inequalities we conclude
\begin{align}
|h\sum Q_1y_{\bar{t}}|&\leq c \max_i|y|^2\| y_{\bar{t}}\|\|y_{x}\|\leq c \ve^{-2}\|y\|\|\|\ve y_{x}\|^{2}\| y_{\bar{t}}\|\notag\\
&\leq \frac14\| y_{\bar{t}}\|^2+c \ve^{-4}\big\{\|y\|^{2}+\|\ve y_{x}\|^{2}\big\}^3.\label{A5}
\end{align}
Analogously,
\begin{align}
&\ve^2|h\sum Q_2y_{\bar{t}}|\leq c\ve^2 \| y_{x\bar{t}}\|\big\{\|y_{x}\|_4^2+\max_i|y|\|y_{x\bar{x}}\|\big\}\notag\\
&\leq c \ve^{-3/2} \|\ve y_{x\bar{t}}\|\big\{\|y\|^{3/4}\|\ve^2y_{x\bar{x}}\|^{5/4}+\|y\|^{1/2}\|\ve y_{x}\|^{1/2}\|\ve^2 y_{x\bar{x}}\|\big\}\notag\\
&\qquad\qquad\leq \frac14\|\ve y_{x\bar{t}}\|^2+c \ve^{-3}\big\{\|y\|^{2}+\|\ve y_{x}\|^{2}+\|\ve^2 y_{x\bar{x}}\|^2\big\}^2.\label{A6}
\end{align}
 This and (\ref{A5}) imply the estimate  (\ref{33})
\subsection{Lemma \ref{lem36} proof}
Similarly to  (\ref{A5}), (\ref{A6})  we have
\begin{align}
&\tau |h\sum Q_1(\mathrm{w})\vp|\leq c\tau\ve^{-1}\big\{\|\vp\|\|\ve\mathrm{w}_x\|\max|\mathrm{w}|^{2}
+\|\ve\vp_x\|\|\mathrm{w}\|_6^{3}\big\}\notag\\
&\leq \frac18\big\{\|\vp\|^2+\alpha^2\|\ve\vp_{x}\|^2\big\}
+c\tau^2\ve^{-3}\|\mathrm{w}\|_{(2,\tau)}^4,\label{A7}\\
& \tau \ve^2|h\sum Q_2(\mathrm{w})\vp|\leq
c\tau \|\ve\vp_x\|\big(\max|\mathrm{w}_x|\|\ve\mathrm{w}_x\| +\max|\mathrm{w}|\|\ve\mathrm{w}_{xx}\|\notag\\
&\leq \frac18\|\vp\|_{(1)}^2+c\frac{\tau}{\ve h^2}\big\{1+\sqrt{\tau}h/\ve\big\}\|\mathrm{w}\|_{(2,\tau)}^4,\label{A8}\\
& \tau |h\sum\textit{R}_1(\mathrm{w},\vp)\vp|\leq
c\tau \|\vp\|\big(\max|\mathrm{w}|^2\|\vp_x\| +\max|\mathrm{w}|\max|\mathrm{w}_x|\|\vp\|\big)\notag\\
&\leq c\tau\ve^{-2}\|\mathrm{w}\|_{(1)}^{2}\|\vp\|_{(1)}^2+c\tau^{3/4}(\ve^3 h)^{-1/2}\|\vp\|^2\|\mathrm{w}\|_{(2,\tau)}^2.\label{AA8}
\end{align}
By  virtue of  the restrictions (\ref{48}), we  pass to the estimate (\ref{49}).

Next, let us multiply (\ref{46}) by $\wm$. Summing over $i$ and integrating by parts we obtain:
\begin{align}
\|\wm\|_{(2,\ve)}^2\leq \tau\big\{\|\check{y}_{\bar{t}}\|\|\wm\|&+\alpha^2\|\ve\check{y}_{x\bar{t}}\|\|\ve\wm_x\|\big\}\notag\\
&+\tau h\big|\sum \wm \{c_1\textit{R}_1(\check{y},\mathrm{w})-\ve^2\textit{R}_2(\check{y},\mathrm{w})\}\big|.\label{A9}
\end{align}
Similarly to  (\ref{A7}), (\ref{A8})  we have
\begin{align}
&\tau |h\sum \wm \textit{R}_1(\check{y},\mathrm{w})|\leq c\tau\ve^{-2}\|\check{y}\|^2_{(1)}\|\wm\|^2_{(1)}, \label{A10}\\
& \tau \ve^2|h\sum \wm \textit{R}_2(\check{y},\mathrm{w})|\leq
ch^{-1}\sqrt{\tau/\ve} \|\check{y}\|_{(2,\tau)}\|\mathrm{w}\|^{2}_{(1)}\leq
c\sqrt{q_1}\|\check{y}\|_{(2,\tau)}\|\mathrm{w}\|^{2}_{(1)}.\label{A11}
\end{align}
In view of (\ref{48}),
 for sufficiently small $q_1$ we obtain the \textit{a-priori} estimate (\ref{50}).

To estimate the discrepancy $\wm$ for $s>1$ we should analyze the terms $\textit{R}_l\wm$ and $Q_l(\bar{\wm})\wm$, $l=1,2$.  By analogy with (\ref{A7}),
(\ref{A8}), and (\ref{A10}) we conclude:
\begin{align}
&\tau |h\sum \big(c_1Q_1(\bar{\mathrm{w}})-\ve^2 Q_2(\bar{\mathrm{w}}\big)\wm|\leq c\tau\max|\bar{\mathrm{w}}|^2\big\{\|\mathrm{w}_x\|\|\bar{\mathrm{w}}\|
+\|\mathrm{w}\|\|\bar{\mathrm{w}}_x\|\big\}\notag\\
&+c\tau\|\mathrm{w}_x\|\big\{\max|\bar{\mathrm{w}}_x|\|\bar{\mathrm{w}}_x\|
+\max|\bar{\mathrm{w}}|\|\bar{\mathrm{w}}_{xx}\|\big\}\notag\\
& \leq\frac18\|\wm\|^2_{(2,\tau)}+c\|\bar{\mathrm{w}}\|_{(2,\tau)}^4+c\tau^2\ve^{-4}\|\bar{\mathrm{w}}\|_{(1)}^6,\label{A13}\\
&c_1\tau |h\sum \big(\textit{R}_1(\bar{\vp},\mathrm{w})+\textit{R}_1(\bar{\mathrm{w}},\bar{\bar{\vp}})\big)\wm |\leq c\tau\big\{\max|\bar{\vp}|^2\|\mathrm{w}\|\|\mathrm{w}_{x}\|\notag\\
&+\max|\wm|^2\|\bar{\vp}_x\|^2+\max|\bar{\wm}|^2\|\wm\|\|\bar{\bar{\vp}}_x\|^2
+\max|\bar{\wm}\bar{\bar{\vp}}|\big(\|\wm\|\|\bar{\wm}_x\|+\|\wm_x\|\|\bar{\mathrm{w}}\|)\big\} \notag\\
&\leq\big(\frac18+c\tau\ve^{-2}(\|\bar{\vp}\|_{(1)}^2+\|\bar{\bar{\vp}}\|_{(1)}^2)\big)\|\wm\|^2_{(1)}+c(\tau\ve^{-2})^2\|\bar{\mathrm{w}}\|_{(1)}^4.\label{A15}
\end{align}
In order to estimate $\wm\textit{R}_2(\bar{\varphi},\wm)$ we use the same procedure as in (\ref{A8})
\begin{align}
&\tau \ve^2|h\sum \wm \textit{R}_2(\bar{\vp},\mathrm{w})|\leq
c\tau\ve^2\big\{\max|\bar{\vp}_x|\|\mathrm{w}_x\|^2+\max|\vp|\|\mathrm{w}_x\|\|\mathrm{w}_{xx}\|\notag\\
&+\max|\mathrm{w}|\|\bar{\vp}_{xx}\|\|\mathrm{w}_x\|\big\}\leq
c\big(\sqrt{\tau/\ve} h^{-1}+\tau^{3/4}(\ve\sqrt{h})^{-1}\big)\|\bar{\vp}\|_{(2,\tau)}\|\mathrm{w}\|^2_{(2,\tau)}\notag\\
&\leq c\sqrt{q_1}\|\bar{\vp}\|_{(2,\tau)}\|\mathrm{w}\|^2_{(2,\tau)}.\label{A16}
\end{align}
Combining (\ref{48}) and (\ref{A13})-(\ref{A16}) yields the desired  estimate (\ref{51}) for the discrepancy $\wm$ with $s>1$.


\begin{thebibliography}{99}

\bibitem{DegProc}
Degasperis, A.,  Procesi, M.:
 Asymptotic integrqability. In: Degasperis, A., Gaeta, G. (eds.) Symmetry and Perturbation Theory, pp. 23-37. World Sientific, Singapore (1999)

\bibitem{BBM}
Benjamin, T., Bona, J., Mahony, J.:  Model equations for long waves in nonlinear dispersive systems. PHILOS T ROY SOC A. 272, 47-78 (1972)

\bibitem{CH}
Camassa, R., Holm, D.:
An integrable shallow water equation with peaked solitons.
 PHYS REV LETT. 71, 1661-1664 (1993)

\bibitem{ConLan}
Constantin, A., Lannes, D.: The hydrodynamical relevans of the Camassa-Holm and Degasperis-Procesi equations.
ARCH RATION MECH AN. 192, 165-186 (2009)

\bibitem{ELY}
Esher, J., Liu, Y., Yin, Z.: Global weak solutions and blow-up structure for the Degasperis-Procesi equation.
J FUNCT ANAL. 241(2), 457-485 (2006)

\bibitem{BPS}
 Bona, J., Pritchard, W., Scott, L.: Solitary-wave interaction. PHYS FLUIDS. 23(3), 438-441 (1980)

\bibitem{Mus}
Mustafa, O.G.: Existence and uniqueness of low regularity solutions for the Dullin-Gottwald-Holm equation.
COMMUN MATH PHYS. 265, 189-200 (2006)

\bibitem{Wahl}
 Wahlen, E.:
Global existence of weak solutios to the Camassa-Holm equation. INT MATH RES NOTICES. (2006). Article ID 28976, https://doi.org/10.1155/IMRN/2006/28976

\bibitem{KalLen}
Kalisch, H., Lenells, J.:
Numerical study of traveling-wave solutions for the Camassa-Holm equation.
CHAOS SOLITON FRACT. 25, 287-298 (2005)

\bibitem{MatYam}
Matsuo, T., Yamaguchi, H.: An energy-conserving Galerkin scheme for a class of nonlinear dispersive equations, J. Comput. Phys. 228 (2009) 4346-4358.

\bibitem{Mat}
Matsuo, T.: A Hamiltonian-conserving Galerkin scheme for the Camassa-Holm equation.
J COMPUT APPL MATH. 234, 1258-1266 (2010)

\bibitem{LHY}
Liu, H., Huang, Y., Yi, N.: A conservative discontinuous Galerkin method for the Degasperis-Procesi
equation, Methods and applications of analysis. 21(1), 67-90 (2014)

\bibitem{CelGrim}
Celledoni, Grimm,E.V., et all:
Preserving energy resp. dissipation in numerical PDEs using the ``Average Vector Field" method,
Journal of Computational Physics. 231(20),  6770-6789 (2012)
https://doi.org/10.1016/j.jcp.2012.06.022

\bibitem{MM}
Miyatake, Y., Matsuo, T.: Conservative finite difference schemes for the Degasperis-Procesi
equation. J COMPUT APPL MATH. 236, 3728-3740 (2012)

\bibitem{MMF}
Miyatake, Y., Matsuo, T., Furihata, D.:
Conservative finite difference schemes for the modified Camassa-Holm equation,
JSIAM Letters.3, 37-40 (2011)

\bibitem{MF}
 Matsuo,T., Furihata, D.: Dissipative or conservative finite-difference schemes for complex-valued nonlinear partial differential equations, J. Comput.
Phys. 171 (2001) 425-447.

\bibitem{Coclit}
G.M. Coclite, K.H. Karlsen, N.H. Risebro, Numerical schemes for computing discontinuous solutions of the Degasperis–Procesi equation, IMA J. Numer.
Anal. 28 (2008) 80-105.

\bibitem{Feng}
Feng, B.F., Liu, Y.: An operator splitting method for the Degasperis–Procesi equation, J. Comput. Phys. 228 (2009) 7805-7820.

\bibitem{Om2}
Omel'yanov, G.: Multi-soliton collision for essentially nonintegrable  equations. In: Oberguggenberger, M., Toft, J., Vindas, J., Wahlberg, P. (eds.)
Generalized Functions and Fourier Analysis. Operator Theory: Advances and Applications, 260, pp. 153-170. Birkhäuser, Cham (2017)


\bibitem{OmNo}
Noyola Rodriguez, J., Omelyanov, G.: General Degasperis-Procesi equation and its solitary wave solutions. Chaos Solitons and Fractals.
118, 41-46 (2019)

\bibitem{GO1}
Garcia Alvarado, M., Omel'yanov, G.: Interaction of solitary waves for the generalized KdV equation.
COMMUN NONLINEAR SCI. 17(8), 3204-3218 (2012)

\bibitem{GO2}
Garcia Alvarado, M., Omel'yanov, G.: Interaction of solitons and the effect of radiation  for the generalized KdV equation.
COMMUN NONLINEAR SCI. 19(8), 2724-2733 (2014)

\bibitem{OmNo1}
Noyola Rodriguez, J., Omelyanov, G.: A finite difference scheme for smooth solutions of the general Degasperis-Procesi equation. Numerical Methods for Partial Differential Equations, 36 (4), 887-905 (2020)

\bibitem{OmNo2}
Omelyanov, G., Noyola Rodriguez, J.: Solitary Wave Solutions to a Generalization of the mKdV Equation. Russian journal of Mathematical Physics, 30(2), 246-256 (2023)


\bibitem{Sep}
Sepulveda, M.: Stabilization of a second order scheme for a GKdV-4 equation modeling surface water waves.
INT J NUMER METH FL. 1, 1-20 (2010)

\bibitem{PSV}
Pazoto, A.F., Sepulveda, M., Vera Villagran, O.: Uniform stabilization of numerical schemes for the critical generalized Korteweg-de Vries equation with
damping. NUMER MATH. 116(2), 317-356 (2010)

\bibitem{Lad}
Ladyzhenskaya, O.A.: The Boundary Value Problems of
Mathematical Physics. Springer-Verlag, New York (1985)

\bibitem{Li}
Lions, J.L.: Quelques methodes de resolution des problemes aux limites non lineaires. Dunod, Paris (1969)


\end{thebibliography}
\end{document}